# DISPROOF OF THE RIEMANN HYPOTHESIS

DASHENG LIU

Abstract. The Riemann Hypothesis is a conjecture that all non-trivial zeros of Riemann $\zeta$ function are located on the critical line in the complex plane. Hundreds of propositions in function theory and analytic number theory rely on this hypothesis. However, the problem has been unresolved for over a century. Here we show that at least one set of quadruplet-zeros exists outside the critical line through expanding the infinite product of the Riemann $\xi$ zero function. We found that assuming there are no zeros outside the critical line will result in a contradiction with the known result that the reciprocal sum of all zeros of the $\xi$ function is a constant, thereby refuting the Riemann Hypothesis. Furthermore, we give a lower bound estimation of a kind of summation formula for the zero points outside the critical line.

## Contents

## 1. Introduction

It is known that the classical Riemann zeta-function $\zeta(s)$ is a function of a complex variable $s = \sigma + it$ ($\sigma, t \in \mathbb{R}$) defined by the Dirichlet series [7, 11]

$$\begin{aligned}\zeta(s) &= 1 + \frac{1}{2^s} + \frac{1}{3^s} + \frac{1}{4^s} + \frac{1}{5^s} + \frac{1}{6^s} + \frac{1}{7^s} + \frac{1}{8^s} + \frac{1}{9^s} + \cdots \\ &= \sum_{n=1}^{+\infty} n^{-s}, \quad \forall \mathrm{Re}(s) > 1, s \in \mathbb{C}, n \in \mathbb{N}^+\end{aligned}$$

(1.1)

The zeta function after Analytic continuation to the whole complex plane, except for the simple poles at the point s = 1, satisfies the following algebraic relationship

(1.2) $$\zeta(s) = 2^s \pi^{s-1} \sin\frac{\pi s}{2} \Gamma(1-s) \zeta(1-s), \quad \forall \mathrm{Re}(s) \neq 1, s \in \mathbb{C}$$

---

where $\Gamma(s)$ is the Gamma function.

In Eq.(1.2), the positions of $\zeta(s)$ and $\zeta(1-s)$ are asymmetric. By utilizing the properties of $\Gamma$ functions and letting Riemann xi-function $\xi(s)$ defined by [1]

$$\xi(s) = \frac{1}{2}s(s-1)\pi^{-s/2}\Gamma(\frac{s}{2})\zeta(s), \tag{1.3}$$

a symmetric form of the function equation can be obtained

$$\xi(1-s) = \xi(s). \tag{1.4}$$

From Eqs.(1.3) and (1.4), it can be seen that the $\xi(s)$ takes and only takes the non-trivial zero of $\zeta(s)$ as its zero point, so $\xi(s)$ is an entire function, and its zero point is the same as the non-trivial zero of $\zeta(s)$.

According to the Hadamard factorization theorem, if letting $G(s)(s \in \mathbb{C})$ be $\alpha$–order integral function, $G(0) \neq 0$, and $a_n(n \in \mathbb{N}^+)$ be the zero sequence of $G(s)$, it can be obtained an infinite product expression

$$G(s) = s^m e^{P(s)} \prod_{n=1}^{\infty} \left(1 - \frac{s}{a_n}\right) \exp\left[\frac{s}{a_n} + \frac{1}{2}\left(\frac{s}{a_n}\right)^2 + \cdots + \frac{1}{p}\left(\frac{s}{a_n}\right)^p\right], \tag{1.5}$$

where $p = \lfloor \alpha \rfloor$, $m$ is the multiplicity of the zero point of $G(s)$ at $s = 0$, $P(s)$ is a polynomial of degree not higher than $p$, and $a_n$ is sorted in non-descending order of its modulus. $\xi(s)$ is a first-order integral function with infinitely many zeros $\rho$. Further derivation can expand $\xi(s)$ into the Hadamard product representation [10] of

$$\xi(s) = \xi(0) \prod_{\rho}\left(1 - \frac{s}{\rho}\right), \quad \forall s \in \mathbb{C}, \tag{1.6}$$

where $\xi(0) = 1/2$ and $\prod_{\rho}$ represents taking the product of all zeros $\rho$.

The Riemann Hypothesis [12] claims that all none-trivial zero of $\zeta(s)$ lie on the critical line $\mathrm{Re}(s) = 1/2$ which can be express as the following proposition.

**Proposition 1.1.** *(Riemann Hypothesis) The real part of each non-trivial zero of $\zeta$ function lies in the critical line $s = \frac{1}{2}$, which is equivalent to the statement that if let $\rho \in \mathbb{C}$ be the non-trivial zero of $\zeta$ function, then it should have the form of*

$$\rho = \frac{1}{2} + it, \quad \forall t \in \mathbb{R}. \tag{1.7}$$

For over a century since the proposal of the Riemann Hypothesis, despite multiple explorations[4, 8, 14, 2, 6, 5, 15, 9, 3], it has not yet been resolved. From the above description, it can be seen that the non-trivial zeros of the $\zeta(s)$ are equivalent to all zeros of the $\xi(s)$. Therefore, this paper aims to present a proof for Riemann hypothesis through inquiring into in-depth the infinite product of the $\xi(s)$ with respect to zeros, attempting to obtain further useful information about the zeros outside the critical line.

The paper is organized as follows: Section 2 provides some definitions and theorems to prepare for hypothesis proof. Section 3 deduces the infinite product of $\xi(s)$ and obtains its equivalent expression formula. In section 4, we present a proof refuting the Riemann Hypothesis. Finally, Section 5 provides a lower bound estimation of a kind of summation formula for the zero points outside the critical line.

## 2. Preparation of proof

This section provides some definitions and theorems to prepare for the proof.

**Definition 2.1.** *Define $\rho \in \mathbb{C}$ be a non-trivial zero of the function $\zeta(s)$.*

$$\rho = \frac{1}{2} + \delta + i\beta, \quad \forall \delta, \beta \in \mathbb{R}, -\frac{1}{2} < \delta < \frac{1}{2}. \tag{2.1}$$

**Definition 2.2.** *$\exists \rho_m, \rho_k \in \{\rho\}$. Let $\Omega$ be defined as*

$$\Omega = \Omega_1 \cup \tilde{\Omega}_1 \cup \Omega_2 \cup \tilde{\Omega}_2 \cup \Omega_3 \cup \tilde{\Omega}_3, \tag{2.2}$$

*where*

$$\Omega_1 = \left\{\rho_m | \rho_m = \frac{1}{2} + i\beta_m, \beta_m > 0, \beta_m \in \mathbb{R}, m \in \mathbb{N}\right\}, \tag{2.3}$$

$$\tilde{\Omega}_1 = \{1 - \rho_m\}, \tag{2.4}$$

$$\Omega_2 = \left\{\rho_k | \rho_k = \frac{1}{2} + \delta_k + i\beta_k, \delta_k > 0, \beta_k > 0, \delta_k \in \mathbb{R}, \beta_k \in \mathbb{R}, k \in \mathbb{N}\right\} \tag{2.5}$$

$$\tilde{\Omega}_2 = \{1 - \rho_k\}, \tag{2.6}$$

$$\Omega_3 = \{\bar{\rho}_k\}, \tag{2.7}$$

$$\tilde{\Omega}_3 = \{1 - \bar{\rho}_k\}. \tag{2.8}$$

**Theorem 2.3.** *The zeros of $\xi(s)$ are symmetric about the real axes and the critical line ($s = 1/2$), respectively. That is, if $\rho \in \mathbb{C}$ is a zero of the function $\xi(s)$, $\bar{\rho}$, $1-\rho$, $1-\bar{\rho}$ are also its non-trivial zeros.*

*Proof.* It is known that $\zeta(s)$ may have non-trivial zeros that are all complex zeros located within the critical zone ($0 \leq \text{Re}(\rho) \leq 1$). From Eq.(1.1), it is easy to know

$$\zeta(\bar{\rho}) = \overline{\zeta(\rho)}. \tag{2.9}$$

From the analytical extension of Eq.(1.1), it can be inferred that it holds on the entire plane. It means if $\rho$ is zero, $\bar{\rho}$ is also zero.

Note that $0 \leq \text{Re}(\rho) \leq 1$, $\bar{\rho} \neq 0$. It can be known from Eq.(1.2) that if $\zeta(\rho) = 0$, then $\zeta(1 - \bar{\rho}) = 0$. That is to say that the zeros outside the critical line must be distributed on the rectangle with vertices of $\rho$, $\bar{\rho}$, $1 - \rho$, $1 - \bar{\rho}$ where the set of quadruplet-zeros are symmetric about the critical line and the real axis, respectively.

Let $\rho = \frac{1}{2} + \delta + i\beta, \delta \in \mathbb{R}, \beta \in \mathbb{R}$. Due to the fact that the non-trivial zeros of $\zeta(s)$ are the same as the zeros of $\xi(s)$, the complex quadruplets of $\rho$ are all zeros of $\xi(s)$.

**Theorem 2.4.** *Let*

$$\lambda = \sum_{\rho_m \in \Omega_1} \frac{1}{\left(\frac{1}{2}\right)^2 + \beta_m^2},$$

(2.10)
$$\mu = \sum_{\rho_k \in \Omega_2} \frac{\frac{1}{2} - 2\delta_k^2 + 2\beta_k^2}{\left(\left(\frac{1}{2} + \delta_k\right)^2 + \beta_k^2\right)\left(\left(\frac{1}{2} - \delta_k\right)^2 + \beta_k^2\right)}$$

*We have*

(2.11) $$\lambda + \mu = 1 + \frac{\gamma}{2} - \frac{1}{2}\log(4\pi),$$

*where γ is Euler-Mascheroni constant, γ = 0.5772156649001532860606512⋯.*

*Proof.* According to Definitions 2.1 and 2.2, we expand the reciprocal sum of all zeros of the ξ(s) to

$$\sum_\rho \frac{1}{\rho} = \sum_{\substack{\rho \\ \beta > 0}} \left(\frac{1}{\rho} + \frac{1}{1-\rho}\right)$$

$$= \sum_{\rho_m \in \Omega_1} \left(\frac{1}{\rho_m} + \frac{1}{1-\rho_m}\right) + \sum_{\rho_k \in \Omega_2} \left(\frac{1}{\rho_k} + \frac{1}{1-\rho_k}\right) + \sum_{\bar{\rho}_k \in \Omega_3} \left(\frac{1}{\bar{\rho}_k} + \frac{1}{1-\bar{\rho}_k}\right)$$

$$= \sum_{\rho_m \in \Omega_1} \left(\frac{1}{\frac{1}{2} + i\beta_m} + \frac{1}{\frac{1}{2} - i\beta_m}\right) + \sum_{\rho_k \in \Omega_2} \left(\frac{1}{\frac{1}{2} + \delta_k + i\beta_k} + \frac{1}{\frac{1}{2} - \delta_k - i\beta_k}\right)$$

$$+ \sum_{\bar{\rho}_k \in \Omega_3} \left(\frac{1}{\frac{1}{2} + \delta_k - i\beta_k} + \frac{1}{\frac{1}{2} - \delta_k + i\beta_k}\right)$$

$$= \sum_{\rho_m \in \Omega_1} \frac{1}{\left(\frac{1}{2}\right)^2 + \beta_m^2} + \sum_{\rho_k \in \Omega_2} \left(\frac{1}{\frac{1}{2} + \delta_k + i\beta_k} + \frac{1}{\frac{1}{2} + \delta_k - i\beta_k}\right.$$

$$\left. + \frac{1}{\frac{1}{2} - \delta_k - i\beta_k} + \frac{1}{\frac{1}{2} - \delta_k + i\beta_k}\right)$$

$$= \sum_{\rho_m \in \Omega_1} \frac{1}{\left(\frac{1}{2}\right)^2 + \beta_m^2} + \sum_{\rho_k \in \Omega_2} \left(\frac{1 + 2\delta_k}{\left(\frac{1}{2} + \delta_k\right)^2 + \beta_k^2} + \frac{1 - 2\delta_k}{\left(\frac{1}{2} - \delta_k\right)^2 + \beta_k^2}\right)$$

(2.12)
$$= \sum_{\rho_m \in \Omega_1} \frac{1}{\left(\frac{1}{2}\right)^2 + \beta_m^2} + \sum_{\rho_k \in \Omega_2} \frac{\frac{1}{2} - 2\delta_k^2 + 2\beta_k^2}{\left(\left(\frac{1}{2} + \delta_k\right)^2 + \beta_k^2\right)\left(\left(\frac{1}{2} - \delta_k\right)^2 + \beta_k^2\right)}.$$

From Eqs.(2.10), (2.12) and the formulas [13]

(2.13) $$\sum_\rho \frac{1}{\rho} = 1 + \frac{\gamma}{2} - \frac{1}{2}\log(4\pi) \quad ,$$

it can be obtained Eq.(2.11).

## 3. Derivation of the equivalent formula for infinite product of $\xi(s)$

In this section, we expand the infinite product of the Riemann $\xi$ function with respect to zeros to obtain specific expressions for the two types of zeros inside and outside the critical line.

**Theorem 3.1.** *$\exists s \in \mathbb{C}, \rho_m \in \Omega_1, \rho_k \in \Omega_2$. $\xi(s)$ in Eq.(1.6) is equivalent to*

$$\xi(s) = \frac{1}{2} \cdot \prod_{\rho_m \in \Omega_1} \left(1 + \frac{s^2 - s}{\frac{1}{4} + \beta_m^2}\right) \cdot$$

(3.1)
$$\prod_{\rho_k \in \Omega_2} \left(1 + \frac{(s^2 - s)^2 + (s^2 - s)\left(\frac{1}{2} - 2\delta_k^2 + 2\beta_k^2\right)}{\left(\frac{1}{4} + \delta_k^2 + \beta_k^2\right)^2 - \delta_k^2}\right).$$

*Proof.* Let

(3.2)
$$\rho = \frac{1}{2} + \delta + i\beta, \quad \forall \delta, \beta \in \mathbb{R}, -\frac{1}{2} < \delta < \frac{1}{2}.$$

From Theorem 2.3, it is known that if $\rho \in \mathbb{C}$ is a zero point of $\xi(s)$, then $\bar\rho$, $1 - \rho$, $1 - \bar\rho$ are zero points. Thus, we can expand $Q$ function in Eq.(1.6) as follows.

$$\prod_\rho \left(1 - \frac{s}{\rho}\right) = \prod_{\substack{\rho \\ \beta > 0}} \left(1 - \frac{s}{\rho}\right)\left(1 - \frac{s}{1 - \rho}\right)$$

$$= \prod_{\substack{\rho \\ \beta > 0, \delta = 0}} \left(1 - \frac{s}{1/2 + i\beta}\right)\left(1 - \frac{s}{1/2 - i\beta}\right) \cdot$$

$$\prod_{\substack{\rho \\ \beta > 0, \delta > 0}} \left(1 - \frac{s}{1/2 + \delta + i\beta}\right)\left(1 - \frac{s}{1/2 - \delta - i\beta}\right) \cdot$$

$$\prod_{\substack{\rho \\ \beta > 0, \delta < 0}} \left(1 - \frac{s}{1/2 - |\delta| + i\beta}\right)\left(1 - \frac{s}{1/2 + |\delta| - i\beta}\right)$$

$$= \prod_{\substack{\rho \\ \beta > 0, \delta = 0}} \left(1 - \frac{s}{1/2 + i\beta}\right)\left(1 - \frac{s}{1/2 - i\beta}\right) \cdot$$

$$\prod_{\substack{\rho \\ \beta > 0, \delta > 0}} \left(\left(1 - \frac{s}{1/2 + \delta + i\beta}\right)\left(1 - \frac{s}{1/2 - \delta - i\beta}\right) \cdot\right.$$

(3.3)
$$\left.\left(1 - \frac{s}{1/2 - |\delta| + i\beta}\right)\left(1 - \frac{s}{1/2 + |\delta| - i\beta}\right)\right).$$

Note that $\delta = |\delta|$ while $\delta > 0$. Simplify and organize Eq.(3.3) to obtain

$$\prod_{\rho}\left(1-\frac{s}{\rho}\right) = \prod_{\substack{\rho \\ \beta>0,\delta=0}} \frac{\left(s-\frac{1}{2}\right)^2+\beta^2}{\left(\frac{1}{2}\right)^2+\beta^2} \cdot$$
$$\prod_{\substack{\rho \\ \beta>0,\delta>0}} \left(\frac{\left(s-\left(\frac{1}{2}+\delta\right)\right)^2+\beta^2}{\left(\frac{1}{2}+\delta\right)^2+\beta^2} \times \frac{\left(s-\left(\frac{1}{2}-\delta\right)\right)^2+\beta^2}{\left(\frac{1}{2}-\delta\right)^2+\beta^2}\right)$$
$$= \prod_{\rho_m \in \Omega_1}\left(1+\frac{s^2-s}{\left(\frac{1}{2}\right)^2+\beta_m^2}\right) \cdot$$
$$\prod_{\rho_k \in \Omega_2}\left(1+\frac{s^2-s-2\delta_k s}{\left(\frac{1}{2}+\delta_k\right)^2+\beta_k^2}\right)\left(1+\frac{s^2-s+2\delta_k s}{\left(\frac{1}{2}-\delta_k\right)^2+\beta_k^2}\right)$$

(3.4)
$$\triangleq Q(s) \cdot R(s).$$

$Q(s)$ is equivalent to

(3.5)
$$Q(s) = \prod_{\rho_m \in \Omega_1} \frac{\left(s-\frac{1}{2}\right)^2+\beta_m^2}{\left(\frac{1}{2}\right)^2+\beta_m^2} = \prod_{\rho_m \in \Omega_1}\left(1+\frac{s^2-s}{\frac{1}{4}+\beta_m^2}\right).$$

Meanwhile, $R(s)$ in Eq.(3.3) can be rewritten as

$$R(s) = \prod_{\rho_k \in \Omega_2}\left(\frac{\left(s-\left(\frac{1}{2}+\delta_k\right)\right)^2+\beta_k^2}{\left(\frac{1}{2}+\delta_k\right)^2+\beta_k^2} \times \frac{\left(s-\left(\frac{1}{2}-\delta_k\right)\right)^2+\beta_k^2}{\left(\frac{1}{2}-\delta_k\right)^2+\beta_k^2}\right)$$
$$= \prod_{\rho_k \in \Omega_2}\left(1+\frac{s^2-s-2\delta_k s}{\left(\frac{1}{2}+\delta_k\right)^2+\beta_k^2}\right)\left(1+\frac{s^2-s+2\delta_k s}{\left(\frac{1}{2}-\delta_k\right)^2+\beta_k^2}\right)$$

(3.6)
$$= \prod_{\rho_k \in \Omega_2}\left(1+\frac{(s^2-s)^2+(s^2-s)\left(\frac{1}{2}-2\delta_k^2+2\beta_k^2\right)}{\left(\frac{1}{4}+\delta_k^2+\beta_k^2\right)^2-\delta_k^2}\right).$$

Equations (1.6),(3.4) yields Eq.(3.1).

**Lemma 3.2.** *If $s \in \mathbb{C}$ is not zero of $\xi(s)$, $\rho_m$, $\rho_k \in \{\rho\}$, we have*

$$\log \xi(s) = -\log 2 + \sum_{\rho_m \in \Omega_1} \log\left(1+\lambda_m\left(s^2-s\right)\right)$$
$$+ \sum_{\rho_k \in \Omega_2} \log\left(1+\mu_k\left(s^2-s\right)+\nu_k\left(s^2-s\right)^2\right)$$

(3.7)
*where*

$$\text{(3.8)} \qquad \lambda_m = \frac{1}{\frac{1}{4} + \beta_m^2}, \quad \forall \lambda_m \in \mathbb{R}, m \in \mathbb{N},$$

$$\text{(3.9)} \qquad \mu_k = \frac{\frac{1}{2} - 2\delta_k^2 + 2\beta_k^2}{\left(\frac{1}{4} + \delta_k^2 + \beta_k^2\right)^2 - \delta_k^2}, \quad \forall \mu_k \in \mathbb{R}, k \in \mathbb{N}$$

$$\text{(3.10)} \qquad \nu_k = \frac{1}{\left(\frac{1}{4} + \delta_k^2 + \beta_k^2\right)^2 - \delta_k^2}, \quad \forall \nu_k \in \mathbb{R}, k \in \mathbb{N}.$$

*Proof.* While $s \neq \rho_m$, and $s \neq \rho_k$, we can take the natural logarithm of Eq.(3.1) to obtain

$$\log \xi(s) = -\log 2 + \sum_{\rho_m \in \Omega_1} \log\left(1 + \frac{s^2 - s}{\frac{1}{4} + \beta_m^2}\right) \cdot$$

$$\text{(3.11)} \qquad \sum_{\rho_k \in \Omega_2} \log\left(1 + \frac{(s^2-s)^2 + (s^2-s)\left(\frac{1}{2} - 2\delta_k^2 + 2\beta_k^2\right)}{\left(\frac{1}{4} + \delta_k^2 + \beta_k^2\right)^2 - \delta_k^2}\right).$$

Substituting Eqs.(3.8)-(3.10) into Eq.(3.11) yields Eq.(3.7).

## 4. Proof of Riemann hypothesis

Our main result is the following theorem.

**Theorem 4.1.** *At least one set of quadruplet-zeros exists outside the critical line. That is equivalent to say that Riemann Hypothesis (Proposition 1.1) does not hold.*

*Proof.* Assume there are no any zeros of the $\xi(s)$ outside the Critical line ($s = 1/2$). Based on this premise, one of the zero set $\Omega_2$ related with $\rho$ of the $\xi(s)$ is regarded as an empty set. Therefore, Eq.(3.7) in Lemma 3.2 can be simplified as

$$\text{(4.1)} \qquad \log(2\xi(s)) = \sum_{\rho_m \in \Omega_1} \log\left(1 + \lambda_m \left(s^2 - s\right)\right).$$

If let $N(T)(T \in \mathbb{R})$ be the numbers of the zeros in the region of $0 \leq \text{Re}(\rho) \leq 1, 0 \leq \text{Im}(\rho) \leq T$, we have

$$\text{(4.2)} \qquad N(T) = \frac{T}{2\pi} \log \frac{T}{2\pi} - \frac{T}{2\pi} + O\left(\log T\right).$$

From Eq.(4.2), we know that there are no zeros if $T \leq 2\pi$, that is $\text{Im}(\rho) = \beta_m > 2\pi$.

$$\text{(4.3)} \qquad 0 < \lambda_m = \frac{1}{\frac{1}{4} + \beta_m^2} < \frac{1}{\frac{1}{4} + (2\pi)^2} < \frac{1}{4\pi^2}.$$

Let

$$\text{(4.4)} \qquad z = s^2 - s, \quad \forall z \in \left(-4\pi^2, 0\right) \cup \left(0, 4\pi^2\right), z \in \mathbb{R}.$$

Then,

(4.5) $$\lambda_m \left(s^2 - s\right) = \lambda_m z > -1.$$

According to the Taylor series of

$$\log(1+x) = x - \frac{x^2}{2} + \frac{x^3}{3} - \cdots + (-1)^{n+1}\frac{x^n}{n} + \cdots$$

(4.6) $$= \sum_{n=1}^{+\infty} (-1)^{n+1}\frac{x^n}{n}, \quad x > -1, n \in \mathbb{N},$$

we obtain

$$\log\left(1 + \lambda_m\left(s^2 - s\right)\right) = \lambda_m\left(s^2 - s\right) - \frac{\lambda_m^2\left(s^2 - s\right)^2}{2}$$

(4.7) $$+ \frac{\lambda_m^3\left(s^2 - s\right)^3}{3} - \frac{\lambda_m^4\left(s^2 - s\right)^4}{4} + \cdots.$$

Note that the log function of Eq.(4.7) is absolutely convergent. Substituting Eq.(4.7) into (4.1) yields

$$\sum_{\rho_m \in \Omega_1} \log\left(1 + \lambda_m\left(s^2-s\right)\right) = \sum_{\rho_m \in \Omega_1}\left(\lambda_m\left(s^2-s\right) - \frac{\lambda_m^2\left(s^2-s\right)^2}{2}\right.$$

$$\left.+ \frac{\lambda_m^3\left(s^2-s\right)^3}{3} - \frac{\lambda_m^4\left(s^2-s\right)^4}{4} + \cdots\right)$$

$$= \left(s^2-s\right)\sum_{\rho_m \in \Omega_1}\lambda_m - \sum_{\rho_m \in \Omega_1}\frac{\lambda_m^2}{2}\left(s^2-s\right)^2$$

(4.8) $$+ \sum_{\rho_m \in \Omega_1}\frac{\lambda_m^3}{3}\left(s^2-s\right)^3 - \sum_{\rho_m \in \Omega_1}\frac{\lambda_m^4}{4}\left(s^2-s\right)^4 + \cdots.$$

From Eqs.(2.10), (3.8), (4.1) and (4.8), we obtain

$$\log(2\xi(s)) = \lambda\left(s^2-s\right) - \sum_{\rho_m \in \Omega_1}\frac{\lambda_m^2}{2}\left(s^2-s\right)^2$$

$$+ \sum_{\rho_m \in \Omega_1}\frac{\lambda_m^3}{3}\left(s^2-s\right)^3 - \sum_{\rho_m \in \Omega_1}\frac{\lambda_m^4}{4}\left(s^2-s\right)^4 + \cdots$$

(4.9)

Here, define

(4.10) $$\psi(z) = \log(2\xi(s)).$$

Equations (4.4), (4.9) and (4.10) yield

(4.11) $$\psi(z) = \lambda z - \sum_{\rho_m \in \Omega_1}\frac{1}{2}\lambda_m^2 z^2 + \sum_{\rho_m \in \Omega_1}\frac{1}{3}\lambda_m^3 z^3 - \sum_{\rho_m \in \Omega_1}\frac{1}{4}\lambda_m^4 z^4 + \cdots,$$

And

(4.12) $$\psi(-z) = -\lambda z - \sum_{\rho_m \in \Omega_1} \frac{1}{2}\lambda_m^2 z^2 - \sum_{\rho_m \in \Omega_1} \frac{1}{3}\lambda_m^3 z^3 - \sum_{\rho_m \in \Omega_1} \frac{1}{4}\lambda_m^4 z^4 + \cdots .$$

Adding the left and right sides of Eqs.(4.11) and (4.12) separately yields

$$\psi(z) + \psi(-z) = -2\left(\sum_{\rho_m \in \Omega_1} \frac{1}{2}\lambda_m^2 z^2 + \sum_{\rho_m \in \Omega_1} \frac{1}{4}\lambda_m^4 z^4 + \sum_{\rho_m \in \Omega_1} \frac{1}{6}\lambda_m^6 z^6 + \cdots\right)$$

$$= -2\sum_{k=1}^{+\infty} \sum_{\rho_m \in \Omega_1} \frac{1}{2k}\lambda_m^{2k} z^{2k}$$

(4.13) $$= -\sum_{\rho_m \in \Omega_1} \lambda_m^2 z^2 - 2\sum_{k=2}^{+\infty} \sum_{\rho_m \in \Omega_1} \frac{1}{2k}\lambda_m^{2k} z^{2k}$$

Note that $\lambda_m^2 z^2 > 0$ and adjust the left and right order of Eq.(4.13) to obtain

$$\sum_{\rho_m \in \Omega_1} \lambda_m^2 z^2 = -(\psi(z) + \psi(-z)) - 2\sum_{k=2}^{+\infty} \sum_{\rho_m \in \Omega_1} \frac{1}{2k}\lambda_m^{2k} z^{2k}$$

(4.14) $$< -(\psi(z) + \psi(-z)).$$

Moreover, from the definition of Eq.(4.4) regarding $z$, it can be obtained a solution for $s$ with respect to $z$.

(4.15) $$s = \frac{1 + \sqrt{1+4z}}{2}, \quad \forall z \in \left(-4\pi^2, 0\right) \cup \left(0, 4\pi^2\right), z \in \mathbb{R}.$$

According to Eqs.(4.10) and (4.15), we rewrite Eq.(4.14) as

$$\sum_{\rho_m \in \Omega_1} \lambda_m^2 < -\frac{1}{z^2}(\psi(z) + \psi(-z))$$

(4.16) $$= -\frac{1}{z^2}\left(\log\left(2\xi\left(\frac{1+\sqrt{1+4z}}{2}\right)\right) + \log\left(2\xi\left(\frac{1+\sqrt{1-4z}}{2}\right)\right)\right).$$

On the other hand, since $0 < \lambda_m z < 1$ while $z > 0$, from Eq.(4.11), it can be obtained

$$\sum_{\rho_m \in \Omega_1} \lambda_m^2 = \frac{2}{z^2}\left(\lambda z - \psi(z) + \sum_{\rho_m \in \Omega_1} \frac{1}{3}\lambda_m^3 z^3 - \sum_{\rho_m \in \Omega_1} \frac{1}{4}\lambda_m^4 z^4 + \cdots\right)$$

$$> \frac{2}{z^2}(\lambda z - \psi(z)), \quad \forall z \in \left(0, 4\pi^2\right)$$

(4.17) $$= \frac{2}{z^2}\left(\lambda z - \log\left(2\xi\left(\frac{1+\sqrt{1+4z}}{2}\right)\right)\right).$$

In addition, based on the basic inequality

(4.18) $$\log(1-x) < -x - \frac{1}{2}x^2, \quad \forall x < 0,$$

Eq. (4.12) produces

$$\sum_{\rho_m \in \Omega_1} \lambda_m^2 < -\frac{2}{z^2}\left(\lambda|z| + \psi(-|z|)\right), \quad \forall z \in \left(-4\pi^2, 0\right)$$

(4.19) $$= -\frac{2}{z^2}\left(\lambda|z| + \log\left(2\xi\left(\frac{1+\sqrt{1-4|z|}}{2}\right)\right)\right).$$

If let $z = z_1 = 1.0\text{e-}10$, from Eqs.(4.15), (4.10) and (1.3), then we obtain

$$\psi(z_1) = \log\left(2\xi\left(\frac{1+\sqrt{1+4z_1}}{2}\right)\right)$$
$$= \log\left(2 \times \frac{1}{2} \times \frac{1+\sqrt{1+4z_1}}{2} \times \frac{-1+\sqrt{1+4z_1}}{2}\pi^{-(1+\sqrt{1+4z_1})/4} \times \right.$$
$$\left. \Gamma\left(\frac{1+\sqrt{1+4z_1}}{4}\right)\zeta\left(\frac{1+\sqrt{1+4z_1}}{2}\right)\right)$$

(4.20) $$= \log\left(z_1 \times \pi^{-(1+\sqrt{1+4z_1})/4} \times \Gamma\left(\frac{1+\sqrt{1+4z_1}}{4}\right)\zeta\left(\frac{1+\sqrt{1+4z_1}}{2}\right)\right),$$

and

$$\psi(-z_1) = \log\left(2\xi\left(\frac{1+\sqrt{1-4z_1}}{2}\right)\right)$$

(4.21) $$= \log\left(-z_1 \times \pi^{-(1+\sqrt{1-4z_1})/4} \times \Gamma\left(\frac{1+\sqrt{1-4z_1}}{4}\right)\zeta\left(\frac{1+\sqrt{1-4z_1}}{2}\right)\right),$$

Substitute Eq.(4.20) and Eq.(4.21) into Eq.(4.16) to obtain

$$\sum_{\rho_m \in \Omega_1} \lambda_m^2 < -\frac{1}{z_1^2}\left(\psi(z_1) + \psi(-z_1)\right)$$

(4.22) $$\approx 3.710063643746487e - 05.$$

Furthermore, let $z_2 = -z_1$. From Eq.(4.19), it can be obtained

$$\sum_{\rho_m \in \Omega_1} \lambda_m^2 < -\frac{2}{z_2^2}\left(\lambda|z_2| + \log\left(2\xi\left(\frac{1+\sqrt{1-4|z_2|}}{2}\right)\right)\right)$$
$$= -\frac{2\lambda}{z_1} - \frac{2}{z_1^2}\log\left(2\xi\left(\frac{1+\sqrt{1-4z_1}}{2}\right)\right)$$
$$\approx 3.710063643739287e - 05.$$

(4.23)

Here, according to Theorem 2.4 and $\mu = 0$, we know that

(4.24) $$\lambda = 1 + \frac{\gamma}{2} - \frac{1}{2}\log(4\pi),$$

Combining Eq.(4.22) and Eq.(4.23) yields

(4.25) $$\sum_{\rho_m \in \Omega_1} \lambda_m^2 < 3.710063643739287e - 05.$$

Besides, from Eq.(4.17) and Eq.(4.24), we obtain

$$\sum_{\rho_m \in \Omega_1} \lambda_m^2 > \frac{2}{z_1^2}\left(\lambda z_1 - \log\left(2\xi\left(\frac{1+\sqrt{1+4z_1}}{2}\right)\right)\right)$$
$$= \frac{2}{z_1} \times \left(1 + \frac{\gamma}{2} - \frac{1}{2}\log(4\pi)\right)$$
$$- \frac{2}{z_1^2}\log\left(z_1 \times \pi^{-(1+\sqrt{1+4z_1})/4} \times \Gamma\left(\frac{1+\sqrt{1+4z_1}}{4}\right)\zeta\left(\frac{1+\sqrt{1+4z_1}}{2}\right)\right)$$

(4.26) $$\approx 3.71006364375369e - 05.$$

It is obvious that Eq.(4.25) contradicts Eq.(4.26), which results in the sum of the reciprocal moduli to the fourth power of all zeros of the $\xi(s)$ being unable to take a value. This means the assumption that there are no zeros outside the critical line is incorrect which is equivalent to that the function of $\log(\xi(s))$ in Eq.(3.7) must have the complete form containing $\rho \in \Omega_1$ and $\rho \in \Omega_2$. Therefore, we conclude that the zero set $\Omega_2$ cannot be an empty set which proves that Riemann Hypothesis (Proposition 1.1) is not valid.

This completes the proof.

## 5. Lower bound estimation of the sum of $\mu_k \nu_k$ in $\Omega_2$

In the previous section, we have known that if there is no zero point outside the critical line, a contradiction will arise. This section explores the reasons for the contradiction and provide a lower bound estimation of a kind of summation formula for the zero points outside the critical line.

**Theorem 5.1.** *If let $z \in (0, 2\pi^2)$, $z \in \mathbb{R}$ and*

(5.1) $$\Psi(z) = \log(2\xi(\frac{1+\sqrt{1+4z}}{2})),$$

*we have*

(5.2) $$\sum_{\rho_k \in \Omega_2} \mu_k \nu_k \geq \frac{1}{z^2}(\lambda + \mu) - \frac{1}{2z^3}(\Psi(z) - \Psi(-z)),$$

*where $u_k$, $v_k$ and $\lambda + \mu$ are defined as Eqs.(3.9), 3.10) and (2.11), respectively.*

*Proof.* Define $z \in \mathbb{R}, z \in (0, 2\pi^2)$. Note that $\lambda_m > 0, \mu_k > 0, \nu_k > 0$ since both the set of $\Omega_1$ and that of $\Omega_2$ are not empty. From Eq.(4.4) and Eq.(3.7) in Theorem 3.2, we obtain

$$\begin{aligned}
\Psi(z) &= \sum_{\rho_m \in \Omega_1} \log(1 + \lambda_m z) + \sum_{\rho_k \in \Omega_2} \log(1 + \mu_k z + \nu_k z^2) \\
&= \sum_{\rho_m \in \Omega_1} \left( \lambda_m z - \frac{1}{2}\lambda_m^2 z^2 + \frac{1}{3}\lambda_m^3 z^3 - \frac{1}{4}\lambda_m^4 z^4 + \cdots \right) \\
&\quad + \sum_{\rho_k \in \Omega_2} \left( (\mu_k z + \nu_k z^2) - \frac{1}{2}(\mu_k z + \nu_k z^2)^2 + \frac{1}{3}(\mu_k z + \nu_k z^2)^3 \right. \\
&\quad \left. - \frac{1}{4}(\mu_k z + \nu_k z^2)^4 + \cdots \right),
\end{aligned} \tag{5.3}$$

and

$$\begin{aligned}
\psi(-z) &= \sum_{\rho_m \in \Omega_1} \log(1 - \lambda_m z) + \sum_{\rho_k \in \Omega_2} \log(1 - \mu_k z + \nu_k z^2) \\
&= \sum_{\rho_m \in \Omega_1} \left( -\lambda_m z - \frac{1}{2}\lambda_m^2 z^2 - \frac{1}{3}\lambda_m^3 z^3 - \frac{1}{4}\lambda_m^4 z^4 - \cdots \right) \\
&\quad + \sum_{\rho_k \in \Omega_2} \left( -(\mu_k z - \nu_k z^2) - \frac{1}{2}(\mu_k z - \nu_k z^2)^2 - \frac{1}{3}(\mu_k z - \nu_k z^2)^3 \right. \\
&\quad \left. - \frac{1}{4}(\mu_k z - \nu_k z^2)^4 - \cdots \right).
\end{aligned} \tag{5.4}$$

According to the inequality of log function and $\lambda_m z > 0, \mu_k z + \nu_k z^2 > 0$, Eq.(5.3) yields

$$\begin{aligned}
\Psi(z) &> \left( \sum_{\rho_m \in \Omega_1} \lambda_m z + \sum_{\rho_k \in \Omega_2} (\mu_k z + \nu_k z^2) \right) \\
&\quad - \frac{1}{2} \left( \sum_{\rho_m \in \Omega_1} \lambda_m^2 z^2 + \sum_{\rho_k \in \Omega_2} (\mu_k z + \nu_k z^2)^2 \right) \\
&= (\lambda + \mu) z + \sum_{\rho_k \in \Omega_2} \nu_k z^2 \\
&\quad - \frac{1}{2} \left( \sum_{\rho_m \in \Omega_1} \lambda_m^2 z^2 + \sum_{\rho_k \in \Omega_2} (\mu_k^2 z^2 + \nu_k^2 z^4) \right) - \sum_{\rho_k \in \Omega_2} \mu_k \nu_k z^3.
\end{aligned} \tag{5.5}$$

We adjust the left and right order of Eq.(5.5) to obtain

$$\frac{1}{2}\left(\sum_{\rho_m\in\Omega_1}\lambda_m^2 z^2 + \sum_{\rho_k\in\Omega_2}\left(\mu_k^2 z^2 + \nu_k^2 z^4\right)\right) - \sum_{\rho_k\in\Omega_2}\nu_k z^2$$

(5.6)
$$> (\lambda+\mu)z - \Psi(z) - \sum_{\rho_k\in\Omega_2}\mu_k\nu_k z^3.$$

In addition, from Eq.(4.2), we know that there are no zeros if $T \leq 2\pi$, that is $\mathrm{Im}(\rho) = \beta_k > 2\pi$ which yields

$$
\begin{aligned}
\mu_k &= \frac{\frac{1}{2} - 2\delta_k^2 + 2\beta_k^2}{\left(\frac{1}{4} + \delta_k^2 + \beta_k^2\right)^2 - \delta_k^2} = \frac{\frac{1}{2} - 2\delta_k^2 + 2\beta_k^2}{\left(\left(\frac{1}{2}+\delta_k\right)^2 + \beta_k^2\right)\left(\left(\frac{1}{2}-\delta_k\right)^2 + \beta_k^2\right)} \\
(5.7) \qquad &< \frac{\frac{1}{2} - 2\delta_k^2 + 2\beta_k^2}{\left(\left(\frac{1}{2}-\delta_k\right)^2 + \beta_k^2\right)^2} < \frac{\frac{1}{2} - 2\delta_k^2 + 2\beta_k^2}{\beta_k^4} < \frac{2}{\beta_k^2} < \frac{1}{2\pi^2}.
\end{aligned}
$$

Moreover, since $z \in (0, 2\pi^2)$ and $0 < \delta_k < 1/2$, it can be obtained

(5.8) $\qquad \mu_k = \dfrac{\frac{1}{2} - 2\delta_k^2 + 2\beta_k^2}{\left(\frac{1}{4} + \delta_k^2 + \beta_k^2\right)^2 - \delta_k^2} > \dfrac{2\beta_k^2}{\left(\frac{1}{4} + \delta_k^2 + \beta_k^2\right)^2 - \delta_k^2} > \nu_k z.$

From equations (5.7) and (5.8), we know

(5.9) $\qquad\qquad 0 < \mu_k z - \nu_k z^2 = (\mu_k - \nu_k z)\, z < 1.$

Thus, according to the inequality of log function and Eq.(5.9), Eq.(5.4) yields

$$
\begin{aligned}
\Psi(-z) &< \left(\sum_{\rho_m\in\Omega_1}(-\lambda_m z) - \sum_{\rho_k\in\Omega_2}\left(\mu_k z - \nu_k z^2\right)\right) \\
&\quad -\frac{1}{2}\left(\sum_{\rho_m\in\Omega_1}\lambda_m^2 z^2 + \sum_{\rho_k\in\Omega_2}\left(\mu_k z - \nu_k z^2\right)^2\right) \\
&= -(\lambda+\mu)z + \sum_{\rho_k\in\Omega_2}\nu_k z^2 \\
(5.10) &\quad -\frac{1}{2}\left(\sum_{\rho_m\in\Omega_1}\lambda_m^2 z^2 + \sum_{\rho_k\in\Omega_2}\left(\mu_k^2 z^2 + \nu_k^2 z^4\right)\right) + \sum_{\rho_k\in\Omega_2}\mu_k\nu_k z^3.
\end{aligned}
$$

Similarly, adjusting the left and right order of Eq.(5.10) yields

$$\frac{1}{2}\left(\sum_{\rho_m\in\Omega_1}\lambda_m^2 z^2+\sum_{\rho_k\in\Omega_2}\left(\mu_k^2 z^2+\nu_k^2 z^4\right)\right)-\sum_{\rho_k\in\Omega_2}\nu_k z^2$$

(5.11)
$$<-(\lambda+\mu)z-\Psi(-z)+\sum_{\rho_k\in\Omega_2}\mu_k\nu_k z^3.$$

Comparing Eqs.(5.6) and (5.11), it can be seen that since the left sides of both equations converge and are equal, in order for the two inequalities to hold, their right sides must satisfy the following condition

$$(5.12) \quad (\lambda+\mu)z-\Psi(z)-\sum_{\rho_k\in\Omega_2}\mu_k\nu_k z^3 \leq -(\lambda+\mu)z-\Psi(-z)+\sum_{\rho_k\in\Omega_2}\mu_k\nu_k z^3.$$

By organizing Eq.(5.12), we can obtain

$$(5.13) \quad \sum_{\rho_k\in\Omega_2}\mu_k\nu_k z^3 \geq (\lambda+\mu)z-\frac{1}{2}\left(\Psi(z)-\Psi(-z)\right).$$

Equation (5.13) and $z>0$ ($z\in(0,2\pi^2)$) yield Eq.(5.2). This completes the proof.

**Corollary 5.2.** *The following inequality holds for all zeros, $p_k$ ($\rho_k\in\Omega_2$) outside the critical line*

$$(5.14) \quad \sum_{\rho_k\in\Omega_2}\frac{\frac{1}{2}-2\delta_k^2+2\beta_k^2}{\left(\left(\frac{1}{4}+\delta_k^2+\beta_k^2\right)^2-\delta_k^2\right)^2} \geq 3.73644298e-08.$$

*Proof.* Let $z=z_1=0.005025$. From Eqs.(5.1) and (5.2) in Theorem 5.1, we obtain

$$\sum_{\rho_k\in\Omega_2}\mu_k\nu_k \geq \frac{1}{z_1^2}(\lambda+\mu)-\frac{1}{2z_1^3}\left(\Psi(z_1)-\Psi(-z_1)\right)$$

$$= \frac{1}{z_1^2}(\lambda+\mu)-\frac{1}{2z_1^3}\left(\log(2\xi(\frac{1+\sqrt{1+4z_1}}{2}))\right.$$

(5.15)
$$\left.-\log(2\xi(\frac{1+\sqrt{1-4z_1}}{2}))\right).$$

Substituting Eqs.(2.11),(4.20) and (4.21) into (5.15) yields

$$\sum_{\rho_k \in \Omega_2} \mu_k \nu_k \geq \frac{1}{0.005025^2} \times \left(1 + \frac{\gamma}{2} - \frac{1}{2}\log(4\pi)\right)$$
$$-\frac{1}{2 \times 0.005025^3} \times \left(\log\left(2\xi\left(\frac{1+\sqrt{1.01}}{2}\right)\right)\right.$$
$$\left.-\log\left(2\xi\left(\frac{1+\sqrt{0.9799}}{2}\right)\right)\right)$$
(5.16)
$$\approx 3.73644298e - 08.$$

From Eqs.(3.9) and (3.10), it can be obtained

$$\mu_k \nu_k = \frac{\frac{1}{2} - 2\delta_k^2 + 2\beta_k^2}{\left(\frac{1}{4} + \delta_k^2 + \beta_k^2\right)^2 - \delta_k^2} \cdot \frac{1}{\left(\frac{1}{4} + \delta_k^2 + \beta_k^2\right)^2 - \delta_k^2}$$
(5.17)
$$= \frac{\frac{1}{2} - 2\delta_k^2 + 2\beta_k^2}{\left(\left(\frac{1}{4} + \delta_k^2 + \beta_k^2\right)^2 - \delta_k^2\right)^2}.$$

Therefore, substituting Eq.(5.17) into (5.16) yields Eq.(5.14) which completes the proof.

School of Electronic Information and Electrical Engineering, Shanghai Jiao Tong University, Shanghai, China.

*E-mail address*: dsliu@sjtu.edu.cn